\magnification=1200
\hsize=11.25cm
\hoffset=1.1cm
\raggedbottom
\frenchspacing
\def\qed{\quad\raise -2pt\hbox{\vrule
\vbox to 10pt{\hrule width 4pt
\vfill\hrule}\vrule}}

\long\def\proclaim #1. #2\endproclaim{\medbreak
{\bf #1.\enspace}{\sl#2}\par\medbreak}

\font\eightrm=cmr8
\font\eighti=cmmi8
\font\eightsy=cmsy8
\font\eightbf=cmbx8
\font\eighttt=cmtt8
\font\eightit=cmti8
\font\eightsl=cmsl8
\font\sixrm=cmr6
\font\sixi=cmmi6
\font\sixsy=cmsy6
\font\sixbf=cmbx6

\skewchar\eighti='177 \skewchar\sixi='177
\skewchar\eightsy='60 \skewchar\sixsy='60

\catcode`\@=11

\def\eightpoint{%
  \textfont0=\eightrm \scriptfont0=\sixrm \scriptscriptfont0=\fiverm
  \def\rm{\fam\z@\eightrm}%
  \textfont1=\eighti \scriptfont1=\sixi \scriptscriptfont1=\fivei
  \def\oldstyle{\fam\@ne\eighti}%
  \textfont2=\eightsy \scriptfont2=\sixsy
\scriptscriptfont2=\fivesy
  \textfont\itfam=\eightit
  \def\it{\fam\itfam\eightit}%
  \textfont\slfam=\eightsl
  \def\sl{\fam\slfam\eightsl}%
  \textfont\bffam=\eightbf \scriptfont\bffam=\sixbf
  \scriptscriptfont\bffam=\fivebf
  \def\bf{\fam\bffam\eightbf}%
  \textfont\ttfam=\eighttt
  \def\tt{\fam\ttfam\eighttt}%
  \abovedisplayskip=9pt plus 2pt minus 6pt
  \abovedisplayshortskip=0pt plus 2pt
  \belowdisplayskip=9pt plus 2pt minus 6pt
  \belowdisplayshortskip=5pt plus 2pt minus 3pt
  \smallskipamount=2pt plus 1pt minus 1pt
  \medskipamount=4pt plus 2pt minus 1pt
  \bigskipamount=9pt plus 3pt minus 3pt
  \normalbaselineskip=9pt
  \setbox\strutbox=\hbox{\vrule height7pt depth2pt width0pt}%
  \let\bigf@ntpc=\eightrm \let\smallf@ntpc=\sixrm
  \normalbaselines\rm}

\def\appeln@te{}
\def\vfootnote#1{\def\@parameter{#1}\insert
  \footins\bgroup\eightpoint
  \interlinepenalty\interfootnotelinepenalty
  \splittopskip\ht\strutbox %
  \splitmaxdepth\dp\strutbox \floatingpenalty\@MM
  \leftskip\z@skip \rightskip\z@skip
  \ifx\appeln@te\@parameter\indent \else{\noindent #1\ }\fi
  \footstrut\futurelet\next\fo@t}

\def\footnoterule{\kern-6\p@
  \hrule width 2truein \kern 5.6\p@} 

\catcode`\@=12

\def\Grille{\setbox13=\vbox to 5mm{\hrule width 110mm\vfill}
\setbox13=\vbox{\offinterlineskip
   \copy13\copy13\copy13\copy13\copy13\copy13\copy13\copy13
   \copy13\copy13\copy13\copy13\box13\hrule width 110mm}
\setbox14=\hbox to 5mm{\vrule height 65mm\hfill}
\setbox14=\hbox{\copy14\copy14\copy14\copy14\copy14\copy14
   \copy14\copy14\copy14\copy14\copy14\copy14\copy14\copy14
   \copy14\copy14\copy14\copy14\copy14\copy14\copy14\copy14\box14}
\ht14=0pt\dp14=0pt\wd14=0pt
\setbox13=\vbox to 0pt{\vss\box13\offinterlineskip\box14}
\wd13=0pt\box13}

\def\fleche(#1,#2)\dir(#3,#4)\long#5{%
\noalign{\nointerlineskip\leftput(#1,#2){\vector(#3,#4){#5}}\nointerlineskip}}

\def\hfl#1#2#3{\smash{\mathop{\hbox to#3{\rightarrowfill}}\limits
^{\scriptstyle#1}_{\scriptstyle#2}}}

\def\gfl#1#2#3{\smash{\mathop{\hbox to#3{\leftarrowfill}}\limits
^{\scriptstyle#1}_{\scriptstyle#2}}}

 \message{`lline' & `vector' macros from LaTeX}
 \catcode`@=11
\def\{{\relax\ifmmode\lbrace\else$\lbrace$\fi}
\def\}{\relax\ifmmode\rbrace\else$\rbrace$\fi}
\def\newcount{\alloc@0\count\countdef\insc@unt}
\def\newdimen{\alloc@1\dimen\dimendef\insc@unt}
\def\newwrite{\alloc@7\write\chardef\sixt@@n}

\newwrite\@unused
\newcount\@tempcnta
\newcount\@tempcntb
\newdimen\@tempdima
\newdimen\@tempdimb
\newbox\@tempboxa

\def\@spaces{\space\space\space\space}
\def\@whilenoop#1{}
\def\@whiledim#1\do #2{\ifdim #1\relax#2\@iwhiledim{#1\relax#2}\fi}
\def\@iwhiledim#1{\ifdim #1\let\@nextwhile=\@iwhiledim
        \else\let\@nextwhile=\@whilenoop\fi\@nextwhile{#1}}
\def\@badlinearg{\@latexerr{Bad \string\line\space or \string\vector
   \space argument}}
\def\@latexerr#1#2{\begingroup
\edef\@tempc{#2}\expandafter\errhelp\expandafter{\@tempc}%
\def\@eha{Your command was ignored.
^^JType \space I <command> <return> \space to replace it
  with another command,^^Jor \space <return> \space to continue without it.}
\def\@ehb{You've lost some text. \space \@ehc}
\def\@ehc{Try typing \space <return>
  \space to proceed.^^JIf that doesn't work, type \space X <return> \space to
  quit.}
\def\@ehd{You're in trouble here.  \space\@ehc}

\typeout{LaTeX error. \space See LaTeX manual for explanation.^^J
 \space\@spaces\@spaces\@spaces Type \space H <return> \space for
 immediate help.}\errmessage{#1}\endgroup}
\def\typeout#1{{\let\protect\string\immediate\write\@unused{#1}}}

\font\tenln    = line10
\font\tenlnw   = linew10

\newdimen\@wholewidth
\newdimen\@halfwidth
\newdimen\unitlength

\unitlength =1pt

\def\thinlines{\let\@linefnt\tenln \let\@circlefnt\tencirc
  \@wholewidth\fontdimen8\tenln \@halfwidth .5\@wholewidth}
\def\thicklines{\let\@linefnt\tenlnw \let\@circlefnt\tencircw
  \@wholewidth\fontdimen8\tenlnw \@halfwidth .5\@wholewidth}

\def\linethickness#1{\@wholewidth #1\relax \@halfwidth .5\@wholewidth}

\newif\if@negarg

\def\lline(#1,#2)#3{\@xarg #1\relax \@yarg #2\relax
\@linelen=#3\unitlength
\ifnum\@xarg =0 \@vline
  \else \ifnum\@yarg =0 \@hline \else \@sline\fi
\fi}

\def\@sline{\ifnum\@xarg< 0 \@negargtrue \@xarg -\@xarg \@yyarg -\@yarg
  \else \@negargfalse \@yyarg \@yarg \fi
\ifnum \@yyarg >0 \@tempcnta\@yyarg \else \@tempcnta -\@yyarg \fi
\ifnum\@tempcnta>6 \@badlinearg\@tempcnta0 \fi
\setbox\@linechar\hbox{\@linefnt\@getlinechar(\@xarg,\@yyarg)}%
\ifnum \@yarg >0 \let\@upordown\raise \@clnht\z@
   \else\let\@upordown\lower \@clnht \ht\@linechar\fi
\@clnwd=\wd\@linechar
\if@negarg \hskip -\wd\@linechar \def\@tempa{\hskip -2\wd\@linechar}\else
     \let\@tempa\relax \fi
\@whiledim \@clnwd <\@linelen \do
  {\@upordown\@clnht\copy\@linechar
   \@tempa
   \advance\@clnht \ht\@linechar
   \advance\@clnwd \wd\@linechar}%
\advance\@clnht -\ht\@linechar
\advance\@clnwd -\wd\@linechar
\@tempdima\@linelen\advance\@tempdima -\@clnwd
\@tempdimb\@tempdima\advance\@tempdimb -\wd\@linechar
\if@negarg \hskip -\@tempdimb \else \hskip \@tempdimb \fi
\multiply\@tempdima \@m
\@tempcnta \@tempdima \@tempdima \wd\@linechar \divide\@tempcnta \@tempdima
\@tempdima \ht\@linechar \multiply\@tempdima \@tempcnta
\divide\@tempdima \@m
\advance\@clnht \@tempdima
\ifdim \@linelen <\wd\@linechar
   \hskip \wd\@linechar
  \else\@upordown\@clnht\copy\@linechar\fi}

\def\@hline{\ifnum \@xarg <0 \hskip -\@linelen \fi
\vrule height \@halfwidth depth \@halfwidth width \@linelen
\ifnum \@xarg <0 \hskip -\@linelen \fi}

\def\@getlinechar(#1,#2){\@tempcnta#1\relax\multiply\@tempcnta 8
\advance\@tempcnta -9 \ifnum #2>0 \advance\@tempcnta #2\relax\else
\advance\@tempcnta -#2\relax\advance\@tempcnta 64 \fi
\char\@tempcnta}

\def\vector(#1,#2)#3{\@xarg #1\relax \@yarg #2\relax
\@linelen=#3\unitlength
\ifnum\@xarg =0 \@vvector
  \else \ifnum\@yarg =0 \@hvector \else \@svector\fi
\fi}

\def\@hvector{\@hline\hbox to 0pt{\@linefnt
\ifnum \@xarg <0 \@getlarrow(1,0)\hss\else
    \hss\@getrarrow(1,0)\fi}}

\def\@vvector{\ifnum \@yarg <0 \@downvector \else \@upvector \fi}

\def\@svector{\@sline
\@tempcnta\@yarg \ifnum\@tempcnta <0 \@tempcnta=-\@tempcnta\fi
\ifnum\@tempcnta <5
  \hskip -\wd\@linechar
  \@upordown\@clnht \hbox{\@linefnt  \if@negarg
  \@getlarrow(\@xarg,\@yyarg) \else \@getrarrow(\@xarg,\@yyarg) \fi}%
\else\@badlinearg\fi}

\def\@getlarrow(#1,#2){\ifnum #2 =\z@ \@tempcnta='33\else
\@tempcnta=#1\relax\multiply\@tempcnta \sixt@@n \advance\@tempcnta
-9 \@tempcntb=#2\relax\multiply\@tempcntb \tw@
\ifnum \@tempcntb >0 \advance\@tempcnta \@tempcntb\relax
\else\advance\@tempcnta -\@tempcntb\advance\@tempcnta 64
\fi\fi\char\@tempcnta}

\def\@getrarrow(#1,#2){\@tempcntb=#2\relax
\ifnum\@tempcntb < 0 \@tempcntb=-\@tempcntb\relax\fi
\ifcase \@tempcntb\relax \@tempcnta='55 \or
\ifnum #1<3 \@tempcnta=#1\relax\multiply\@tempcnta
24 \advance\@tempcnta -6 \else \ifnum #1=3 \@tempcnta=49
\else\@tempcnta=58 \fi\fi\or
\ifnum #1<3 \@tempcnta=#1\relax\multiply\@tempcnta
24 \advance\@tempcnta -3 \else \@tempcnta=51\fi\or
\@tempcnta=#1\relax\multiply\@tempcnta
\sixt@@n \advance\@tempcnta -\tw@ \else
\@tempcnta=#1\relax\multiply\@tempcnta
\sixt@@n \advance\@tempcnta 7 \fi\ifnum #2<0 \advance\@tempcnta 64 \fi
\char\@tempcnta}

\def\@vline{\ifnum \@yarg <0 \@downline \else \@upline\fi}

\def\@upline{\hbox to \z@{\hskip -\@halfwidth \vrule
  width \@wholewidth height \@linelen depth \z@\hss}}

\def\@downline{\hbox to \z@{\hskip -\@halfwidth \vrule
  width \@wholewidth height \z@ depth \@linelen \hss}}

\def\@upvector{\@upline\setbox\@tempboxa\hbox{\@linefnt\char'66}\raise
     \@linelen \hbox to\z@{\lower \ht\@tempboxa\box\@tempboxa\hss}}

\def\@downvector{\@downline\lower \@linelen
      \hbox to \z@{\@linefnt\char'77\hss}}

\thinlines

\newcount\@xarg
\newcount\@yarg
\newcount\@yyarg
\newcount\@multicnt
\newdimen\@xdim
\newdimen\@ydim
\newbox\@linechar
\newdimen\@linelen
\newdimen\@clnwd
\newdimen\@clnht
\newdimen\@dashdim
\newbox\@dashbox
\newcount\@dashcnt
 \catcode`@=12

\newbox\tbox
\newbox\tboxa

\def\leftzer#1{\setbox\tbox=\hbox to 0pt{#1\hss}%
     \ht\tbox=0pt \dp\tbox=0pt \box\tbox}

\def\rightzer#1{\setbox\tbox=\hbox to 0pt{\hss #1}%
     \ht\tbox=0pt \dp\tbox=0pt \box\tbox}

\def\centerzer#1{\setbox\tbox=\hbox to 0pt{\hss #1\hss}%
     \ht\tbox=0pt \dp\tbox=0pt \box\tbox}

\def\image(#1,#2)#3{\vbox to #1{\offinterlineskip
    \vss #3 \vskip #2}}

\def\leftput(#1,#2)#3{\setbox\tboxa=\hbox{%
    \kern #1\unitlength
    \raise #2\unitlength\hbox{\leftzer{#3}}}%
    \ht\tboxa=0pt \wd\tboxa=0pt \dp\tboxa=0pt\box\tboxa}

\def\rightput(#1,#2)#3{\setbox\tboxa=\hbox{%
    \kern #1\unitlength
    \raise #2\unitlength\hbox{\rightzer{#3}}}%
    \ht\tboxa=0pt \wd\tboxa=0pt \dp\tboxa=0pt\box\tboxa}

\def\centerput(#1,#2)#3{\setbox\tboxa=\hbox{%
    \kern #1\unitlength
    \raise #2\unitlength\hbox{\centerzer{#3}}}%
    \ht\tboxa=0pt \wd\tboxa=0pt \dp\tboxa=0pt\box\tboxa}

\unitlength=1mm

\def\cput(#1,#2)#3{\noalign{\nointerlineskip\centerput(#1,#2){#3}
                             \nointerlineskip}}

\def\segment(#1,#2)\dir(#3,#4)\long#5{%
\leftput(#1,#2){\lline(#3,#4){#5}}}

\null
\bigskip

\centerline{\bf A classic proof of a recurrence}
\smallskip
\centerline{\bf for a very classical sequence}
\bigskip
\centerline{Dominique
Foata\/\footnote{$^1$}{D\'epartement de math\'ematique,
Universit\'e Louis Pasteur, 7, rue Ren\'e-Descartes,\hfil\break
F-67084 Strasbourg, France. {\eighttt
(foata@math.u-strasbg.fr)}} and Doron
Zeilberger\/\footnote{$^2$}{Department of mathematics, Temple
University, Philadelphia, Pennsylvania 19122.\hfil\break {\eighttt
(zeilberg@math.temple.edu)}}}

\bigskip
\rightline{\sl To Marco Sch\"utzenberger, in memoriam.\qquad}
\bigskip\bigskip

Richard Stanley [St96] has recently narrated the fascinating
story of how the {\it classical} Schr\"oder [Sch1870] numbers $s(n)$ are
even more classical than has been believed before. They
(at least $s(10)=103049$) have been known to Hipparchus (190-127
{\eightrm B.C.}). Stanley recalled the three-term linear recurrence
([Co64]; [Co74], p. 57)
 $$\displaylines{(1)\qquad
3(2n-1)s(n)=(n+1)s(n+1)+(n-2)s(n-1)\quad (n\ge 2),\hfill\cr
s(1)=s(2)=1,\cr}
$$
and stated that ``no direct combinatorial
proof of this formula seems to be known." The purpose of this
note is to fill this gap.

The present proof reflects the ideas
of our great master, Marcel-Paul Sch\"utzenberger (1920-1996), who
taught us that every algebraic relation is to be given a combinatorial
counterpart and vice versa. The methodology has been vigorously and
successfully pursued by the \'Ecole bordelaise (e.g., [Cor75],
[Vi85].)

The recurrence (1) is tantalizingly similar to the linear recurrence
$$
2(2n-1)c(n)=(n+1)c(n+1)\quad (n\ge 1),\qquad c(1)=1,\leqno(2)
$$
that is obviously satisfied by the Catalan numbers
$c(n)=(1/n){2n-2\choose n-1}$.
Our proof is inspired by R\'emy's elegant combinatorial proof [Re85]
of~(2) shown to us by Viennot [Vi82].

Recall ([Co74], pp. 56-57) that a {\it Schr\"oder tree} $T$ is either
the tree consisting of its root alone $T=r$, or an ordered tuple
$[r;T_1, \ldots , T_l]$, where $l \geq 2$ and $T_1$, \dots~, $T_l$
are smaller Schr\"oder trees. The first symbol~$r$ is called the
{\it root} of~$T$ and the roots of $T_1$, \dots~, $T_l$ are called the
{\it sons} of~$T$. A son-less node is called a {\it leaf}. The number of
Schr\"oder trees with $n$ leaves is denoted by $s(n)$.

The first values of the $s(n)$'s appear in Table~1.
$$\vbox{\halign{\strut\hfil$#$\hfil\quad \vrule\quad
&\hfil$#$\hfil\quad &\hfil$#$\hfil\quad &\hfil$#$\hfil\quad
&\hfil$#$\hfil\quad &\hfil$#$\hfil\quad &\hfil$#$\hfil\quad
&\hfil$#$\hfil\quad &\hfil$#$\hfil\quad
&\hfil$#$\hfil\quad &\hfil$#$\hfil\cr
n&1&2&3&4&5&6&7&8&9&10\cr
\noalign{\hrule}
s(n)&1&1&3&11&45&197&903&4279&20793&103049\cr
\noalign{\medskip}
\multispan{11}\hfil Table 1\hfil\cr}}
$$

\goodbreak
Our proof will be based on another combinatorial model for the
Schr\"oder numbers, the {\it well-weighted binary plane
trees}. Recall that a plane tree is said to be {\it binary}, if each node
has either no sons or exactly two sons. As is well-known
(see, e.g., [Co74], pp.~52-53), the number
of such trees with exactly~$n$ leaves is the {\it Catalan number}
$c(n)$. A node is said to be {\it interior}, if it is not a leaf. The binary
trees we will be using are {\it weighted}. This means that each {\it
interior} node is given a weight equal to~1 or equal to~2. A node is
said to be {\it well-weighted}, if, whenever it has weight~2, its right
son is not a leaf.

A binary plane tree is {\it well-weighted} if it is
weighted and if all its interior nodes are well-weighted. In
short, we will speak of a {\it well-weighted tree}. As shown in
Table~2, there are $c(4)=5$ binary plane trees with four leaves, and
$s(4)=11$ well-weighted trees with four leaves. The nodes that can get
{\it either} weight~1 or weight~2 are indicated by the symbol
``$\bigcirc$."

\centerline{\vbox{\vskip 18mm
\offinterlineskip
\segment(0,0)\dir(-1,1)\long{10}
\segment(0,0)\dir(1,1)\long{10}
\segment(-5,5)\dir(1,2)\long{2.5}
\segment(5,5)\dir(-1,2)\long{2.5}
\leftput(-7,2){1}
\rightput(7,2){1}
\leftput(-1.5,-.5){$\bigcirc$}
\vskip10mm
}\hskip3cm
\vbox{\vskip 18mm
\offinterlineskip
\segment(0,0)\dir(-1,1)\long{15}
\segment(0,0)\dir(1,1)\long{5}
\segment(-5,5)\dir(1,1)\long{5}
\segment(-10,10)\dir(1,1)\long{5}
\rightput(-12,8){1}
\rightput(-7,3){1}
\rightput(0,-3){1}
\vskip10mm
}\hskip3cm
\vbox{\vskip 18mm
\offinterlineskip
\segment(0,0)\dir(-1,1)\long{10}
\segment(0,0)\dir(1,1)\long{5}
\segment(-5,5)\dir(1,1)\long{10}
\segment(0,10)\dir(-1,1)\long{5}
\rightput(1,6){1}
\rightput(0,-3){1}
\leftput(-6.5,4.5){$\bigcirc$}
\vskip10mm
}}

\centerline{\vbox{\vskip 10mm
\offinterlineskip
\segment(0,0)\dir(-1,1)\long{5}
\segment(0,0)\dir(1,1)\long{10}
\segment(5,5)\dir(-1,1)\long{10}
\segment(0,10)\dir(1,1)\long{5}
\rightput(-1,7){1}
\leftput(-1.5,-.5){$\bigcirc$}
\rightput(7,2){1}
\vskip5mm
}\hskip3cm
\vbox{\vskip 10mm
\offinterlineskip
\segment(0,0)\dir(-1,1)\long{5}
\segment(0,0)\dir(1,1)\long{15}
\segment(5,5)\dir(-1,1)\long{5}
\segment(10,10)\dir(-1,1)\long{5}
\leftput(-1.5,-.5){$\bigcirc$}
\leftput(3.5,4.5){$\bigcirc$}
\rightput(12,7){1}
\vskip4mm
}\hskip1cm}

\centerline{Table 2}

\smallskip
A well-weighted tree can be also defined recursively as follows. It is
either a single unweighted node (serving both as root and leaf),
or a triple $[r;T_1, T_2]$, where $r$ is either $1$ or $2$, and
where $T_1$ and $T_2$ are smaller well-weighted trees, with the
provision that if $T_2$ is a mere leaf, then $r$ must~be~$1$.

Define a mapping from the set of Schr\"oder trees to the set of
well-weighted trees as follows: if $T=r$, then
$\Phi(T):=r$; if the root of ~$T$ has exactly two sons, i.e.,
$T=[r;T_1,T_2]$, then $\Phi(T):=[1;\Phi(T_1),\Phi(T_2)]$;
if the root of~$T$ has more than two sons, i.e.,
$T=[r;T_1,\ldots ,T_l]$ with $l>2$, then
$\Phi(T):=[2;\Phi(T_1),\Phi([r; T_2, \ldots , T_l])\,]$.
It is clear that $\Phi$ is a bijection that preserves the number of
leaves.

Being binary, every well-weighted tree with $n$ leaves, possesses
exactly $(n-1)$ interior nodes and hence $(2n-1)$ nodes altogether. A
well-weighted tree with $n$ leaves is said to be {\it pointed}, if
exactly one of its $(2n-1)$ nodes is pointed; it is {\it leaf-pointed}
if exactly one of its~$n$ leaves is pointed; it is {\it interior-pointed}
if exactly one of its $(n-1)$ interior nodes is pointed. Let $PT(n)$,
(resp. $LT(n)$, resp. $IT(n)$) be the set of pointed (resp. leaf-pointed,
resp. interior-pointed) well-weighted trees with $n$ leaves.

\goodbreak
As $(2n-1)s(n)=\# PT(n)$,  $(n+1)s(n+1)=\# LT(n+1)$ and
$(n-2)s(n-1)=\# IT(n-1)$, formula (1) will be proved
combinatorially (or bijectively, as some people say to-day), if we
can construct a bijection~$\sigma$ of $\{1,2,3\}\times PT(n)$
onto the disjoint union $LT(n+1)\cup IT(n-1)$.

\medskip
The construction of such a bijection will consist of adding a new leaf
to each pointed well-weighted tree of $PT(n)$, in three different
ways denoted by $L_1$, $L_2$, $R_1$. We shall get all of $LT(n+1)$
plus a set $B(n+1)$ of leaf-pointed weighted trees, but not
well-weighted, which is in one-to-one correspondence with $IT(n-1)$.

\medskip
To construct the bijection $\sigma$ we proceed as follows. Start
with $t$ in $PT(n)$ and let $s$ be the subtree of~$t$ whose root is the
pointed node of~$t$. For $j=1,2$ define $\sigma'(L_j,t)$ to be the
leaf-pointed weighted tree with $(n+1)$ leaves obtained from~$t$, by
performing the following replacement:

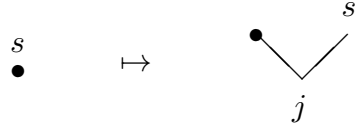
\centerline{\vbox{\vskip 13mm
\offinterlineskip
\leftput(0,0){$\bullet$}
\leftput(0,3){$s$}
\vskip8mm}
\hskip 1cm
\vbox{\vskip 16mm
\offinterlineskip
\leftput(1,4){$\mapsto$}
\vskip5mm}
\hskip 2cm
\vbox{\vskip 13mm
\offinterlineskip
\segment(0,0)\dir(-1,1)\long{5}
\segment(0,0)\dir(1,1)\long{5}
\leftput(-6,4){$\bullet$}
\rightput(6,7){$s$}
\rightput(0.5,-4){$j$}
\vskip8mm
}\kern15mm}

\noindent
Notice that the new interior node receives weight $j$ and the point
gets moved from the root of~$s$ to the new leaf.

In the same manner $\sigma'(R_1,t)$ is defined by performing the
replacement:

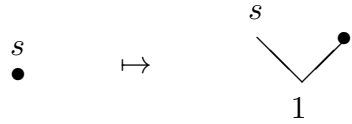
\centerline{\vbox{\vskip 13mm
\offinterlineskip
\leftput(0,0){$\bullet$}
\leftput(0,3){$s$}
\vskip8mm}
\hskip 1cm
\vbox{\vskip 16mm
\offinterlineskip
\leftput(1,4){$\mapsto$}
\vskip5mm}
\hskip 2cm
\vbox{\vskip 13mm
\offinterlineskip
\segment(0,0)\dir(-1,1)\long{5}
\segment(0,0)\dir(1,1)\long{5}
\leftput(-6,7){$s$}
\rightput(5.5,4){$\bullet$}
\rightput(0.5,-4){$1$}
\vskip8mm
}\kern15mm}

\smallskip
Clearly, $\sigma'$ is an {\it injection} of
$\{L_1,L_2,R_1\}\times PT(n)$
into the set of leaf-pointed weighted trees with $(n+1)$ leaves.
Furthermore,
$\sigma'(L_1,t)$ and $\sigma'(R_1,t)$ are always
{\it well-weighted} and  $\sigma'(L_2,t)$ is well-weighted when
the subtree $s$ is {\it not} a leaf. Define $\sigma=\sigma'$ in all
those cases.

When $s$ is a leaf, the subtree of~$t$ whose root is the {\it father} of
the pointed node is one of the following forms

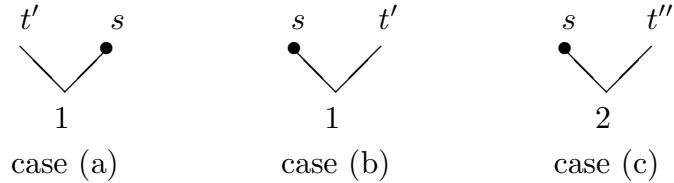
\centerline{\vbox{\vskip 16mm
\offinterlineskip
\segment(0,0)\dir(-1,1)\long{5}
\segment(0,0)\dir(1,1)\long{5}
\leftput(-5,7){$t'$}
\rightput(5.5,4){$\bullet$}
\rightput(0.5,-4){$1$}
\leftput(5,7){$s$}
\centerput(0,-9){case (a)}
\vskip13mm
}\hskip3cm
\vbox{\vskip 16mm
\offinterlineskip
\segment(0,0)\dir(-1,1)\long{5}
\segment(0,0)\dir(1,1)\long{5}
\leftput(-5.5,4){$\bullet$}
\rightput(7,7){$t'$}
\leftput(-5,7){$s$}
\rightput(0.5,-4){$1$}
\centerput(0,-9){case (b)}
\vskip13mm
}\hskip3cm
\vbox{\vskip 16mm
\offinterlineskip
\segment(0,0)\dir(-1,1)\long{5}
\segment(0,0)\dir(1,1)\long{5}
\leftput(-5.5,4){$\bullet$}
\rightput(7.5,7){$t''$}
\leftput(-5,7){$s$}
\rightput(0.5,-4){$2$}
\centerput(0,-9){case (c)}
\vskip13mm
}
}

\noindent
where the subtree $t'$ may be any tree, while the subtree $t''$ must
{\it not} be a leaf. When applying $\sigma'$ to $(L_2,t)$, we will get

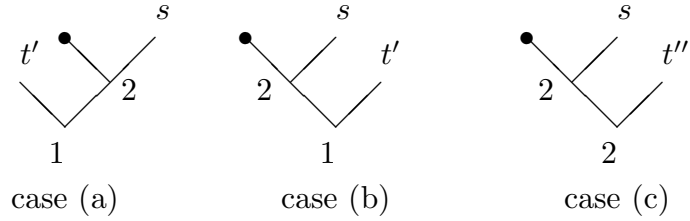
\centerline{\vbox{\vskip 15mm
\offinterlineskip
\segment(0,0)\dir(-1,1)\long{5}
\segment(0,0)\dir(1,1)\long{10}
\segment(5,5)\dir(-1,1)\long{5}
\leftput(-5,7){$t'$}
\rightput(8,3){$2$}
\leftput(10,12){$s$}
\rightput(0,-4){$1$}
\centerput(0,9){$\bullet$}
\centerput(0,-9){case (a)}
\vskip12mm
}\hskip3cm
\vbox{\vskip 15mm
\offinterlineskip
\segment(0,0)\dir(-1,1)\long{10}
\segment(0,0)\dir(1,1)\long{5}
\segment(-5,5)\dir(1,1)\long{5}
\leftput(5,7){$t'$}
\rightput(-7,3){$2$}
\rightput(0,-4){$1$}
\centerput(-10,9){$\bullet$}
\leftput(0,12){$s$}
\centerput(0,-9){case (b)}
\vskip12mm
}
\hskip3cm
\vbox{\vskip 15mm
\offinterlineskip
\segment(0,0)\dir(-1,1)\long{10}
\segment(0,0)\dir(1,1)\long{5}
\segment(-5,5)\dir(1,1)\long{5}
\leftput(5,7){$t''$}
\rightput(-7,3){$2$}
\rightput(0,-4){$2$}
\centerput(-10,9){$\bullet$}
\leftput(0,12){$s$}
\centerput(0,-9){case (c)}
\vskip12mm
}}

\noindent
Alas, none of those trees is well-weighted (since $s$ is a leaf). To
remedy cases (a) and (b) replace them by:

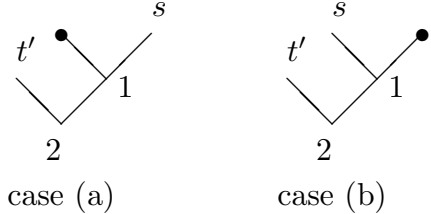
\centerline{\vbox{\vskip 15mm
\offinterlineskip
\segment(0,0)\dir(-1,1)\long{5}
\segment(0,0)\dir(1,1)\long{10}
\segment(5,5)\dir(-1,1)\long{5}
\leftput(-5,7){$t'$}
\rightput(8,3){$1$}
\rightput(0,-4){$2$}
\centerput(0,9){$\bullet$}
\leftput(10,12){$s$}
\centerput(0,-9){case (a)}
\vskip12mm
}\hskip3cm
\vbox{\vskip 15mm
\offinterlineskip
\segment(0,0)\dir(-1,1)\long{5}
\segment(0,0)\dir(1,1)\long{10}
\segment(5,5)\dir(-1,1)\long{5}
\leftput(-5,7){$t'$}
\rightput(8,3){$1$}
\rightput(0,-4){$2$}
\centerput(10,9){$\bullet$}
\leftput(0,12){$s$}
\centerput(0,-9){case (b)}
\vskip12mm
}
}

\noindent
and define $\sigma(L_2,t)$ to be the tree thereby obtained.
It now belongs to $LT(n+1)$. It is
straightforward to verify that $\sigma$ is a {\it bijection} of
the set of all pairs in $\{L_1,L_2,R_1\}\times PT(n)$ that do not belong
to case~(c) onto $LT(n+1)$.

Finally, to obtain $\sigma(L_2,t)$ in case (c), replace the portion
depicted above that belongs to $\sigma'(L_2,t)$
 by just the subtree $t''$ (losing two leaves) and make the root of
$t''$ be the pointed node of the new tree. We get an
interior-pointed well-weighted plane tree with $(n-1)$ leaves.
Furthermore, the restriction of~$\sigma$ to case~(c) is a bijection
onto $IT(n-1)$.

\vglue 38pt
\centerline{\bf References}
{\eightpoint

\bigskip
\def\pointir{\discretionary{.}{}{.\kern1em}\nobreak
\hskip 0em plus .3em minus .4em }

\def\article#1|#2|#3|#4|#5|#6|#7|
    {{\leftskip=7mm\noindent
     \hangindent=2mm\hangafter=1
     \llap{[#1]\hskip.35em}{#2}\pointir
     #3, {\sl #4}, vol.\nobreak\ {\bf #5}, {\oldstyle #6},
     p.\nobreak\ #7.\par}}
\def\livre#1|#2|#3|#4|
    {{\leftskip=7mm\noindent
    \hangindent=2mm\hangafter=1
    \llap{[#1]\hskip.35em}{#2}\pointir
    {\sl #3}\pointir   #4.\par}}
\def\divers#1|#2|#3|
    {{\leftskip=7mm\noindent
    \hangindent=2mm\hangafter=1
     \llap{[#1]\hskip.35em}{#2}\pointir
     #3.\par}}

\article Co64|Louis Comtet|Calcul pratique des coefficients
de Taylor d'une fonction alg\'e\-brique|Enseignement
Math.|10|1964|267--270|

\livre Co74|Louis Comtet|Advanced
Combinatorics|Dordrecht-Holland/Boston, {\oldstyle 1974}|

\divers Cor75|Robert Cori|Un code pour les graphes planaires et ses
applications, {\sl Ast\'erisque}, vol.~{\bf 27}, {\oldstyle 1975}|

\article Re85|Jean-Luc R\'emy|Un proc\'ed\'e it\'eratif de
d\'enombrement d'arbres binaires et son application \`a leur
g\'en\'eration al\'eatoire|RAIRO Inform. Th\'eor.|19|1985|179--195|

\article Sch1870|Ernst Schr\"oder|Vier combinatorische Probleme|Z. f\"ur
Math. Physik|15|1870|361-370|

\divers St96|Richard P. Stanley|Hipparchus,
Plutarch, Schr\"oder and Hough, to appear in {\sl Amer. Math. Monthly},
{\oldstyle 1996}, also available from:\hfil\break
{\tt http://www-math.mit.edu/$\sim$rstan/papers.html}|

\divers Vi82|G\'erard-Xavier Viennot|Oral communication, Oberwolfach,
{\oldstyle 1982}|

\divers Vi85|G\'erard-Xavier Viennot|Probl\`emes combinatoires pos\'es
par la physique statistique, S\'eminaire Bourbaki, 36i\`eme ann\'ee,
{\oldstyle 1983}/{\oldstyle 1984}, no. 626, {\sl Ast\'erisque},
vol. {\bf 225-226}, {\oldstyle 1985}|

}
\bye